\documentclass[12pt]{article}
\usepackage{epsf}
\usepackage{times}
\begin{document}

\newcommand{\noi}{\noindent}

\title{An Experiment with Hierarchical Bayesian Record Linkage} 
\author{Michael D.\ Larsen \\
 The George Washington University \\
 Department of Statistics \\
  mlarsen@bsc.gwu.edu}
\date{December 20, 2012}

\maketitle

\begin{abstract}
In record linkage (RL), or exact file matching, the goal is to identify the
links between entities with information on two or more files. RL
is an important activity in areas including counting the
population, enhancing survey frames and data, and conducting
epidemiological and follow-up studies.  RL 
is challenging when files are very large, no accurate personal
identification (ID) number is present on all files for all units, and some
information is recorded with error.  Without an unique ID number one
must rely on comparisons of names, addresses, dates, and other
information to find the links.  Latent class models can be used to
automatically score the value of information for determining match
status.  Data for fitting models come from comparisons made within
groups of units that pass initial file blocking requirements.  Data
distributions can vary across blocks.  This article examines the use
of prior information and hierarchical latent class models in the
context of RL.   
\end{abstract}


\section{Introduction}
Record linkage (RL; Fellegi and Sunter 1969)  is the name given
to the activity of using information in two or more data bases to
identify units or individuals represented in more than one of the
files. One performs RL when merging files to avoid
duplicate records and to correctly associate information present on
the two or more files with unique individuals.  RL can be
used as part of population size estimation, survey frame and survey
data enhancement, and epidemiological and longitudinal studies. 
In general when there are multiple files on a single population one
could consider using RL for a number of activities. 
Larsen (2012) reviews some literature on the subject.   See also
Winkler (1995), Alvey and Jamerson (1997),  and Herzog, Scheuren, and
Winkler (2007).  

When files are large, accurate and unique personal
identification (ID) numbers are not present on all files for all
units, and there are errors in some information RL
can be non trivial.   Without an unique ID number one
must compare names, addresses, dates, and other
information to find the links between files.  The outcome of
comparisons between two records, one on each of two files, is called a
vector of comparisons, or comparison vector.   Typcially comparisons
are only made for records that pass initial file blocking
requirements.  The pairs that fall outside a block are assumed to be
nonmatches. 

Once a comparison vector is computed one must decide how much evidence
the set of comparisons give in favor of the pair being a match.  In
some applications, predetermined scoring procedures are used.  As an
alternative, one
could imagine taking a sample of cases, determining whether the pairs
truly are matches, and fitting a statistical model to the
results. Such an approach is a type of supervised learning similar in
nature to discriminant analysis or logistic regression.   Latent class
models, which are a type of mixture model, can be fit to used to
automatically score the value of information for determining match
status.   Latent class analysis (LCA) clusters the data and is a form of
unsupervised learning.  Data for fitting models come from comparisons
made within  groups of units that pass initial file blocking
requirements.   LCA has been used in RL by Larsen (2012), Lahiri and
Larsen (2005), Larsen and Rubin (2001), Winkler (1988,  1994), 
Jaro (1989, 1995), and others.  

This article examines the use of prior information and hierarchical
latent class models in the context of RL.  First, the impact of an
informative prior distribution developed from record linkage
operations similar to the current one is studied.  Second, since data
distributions can vary across blocks, a hierarchical latent class
model is used to account for inter-block heterogeneity.  Both
developments are studied through simulation and compared to LCA in
terms of RL error rates. 

Section~\ref{sect:RL} presents statistical models for record linkage.
Section~\ref{sect:comp} discusses computational algorithms. 
Section~\ref{sect:sim} reports on a simulation study. 
Section~\ref{sect:concl} gives a summary and conclusions.

\section{Record Linkage Latent Class Models
 \label{sect:RL}}

Suppose that there are two files, $A$ and $B$, on a single
population.  Consider record $a$ in file $A$ and  record $b$ in $B$.
Do records $a$ and $b$ correspond to the same person or entity?  
Assume files $A$ and $B$ do not contain unique
identification numbers for any units in the files.  Variables in the
two files are used to judge the similarity of the record pairs.  To do
so one defines agreement for each piece of information common to both
files.   In a household-based study, variables can include 
last name currently and at birth, middle name or initial, first name,
house and unit number,  street name, 
age or date of birth, sex, race/ethnicity, and relation to head of
household.   Files often are preprocessed before linkage is attempted.
For example,  names can be standardized and coded according to Soundex
codes or other scheme.   Names and address fields are parsed and
standardized.  Birth date can be separated into day, month, and year. 

In the case of simple comparisons, for each pair of records $(a,b)$
being considered, a vector of 1's and 0's indicating agreement and
disagreement on $K$ comparison fields is recorded.   That is, 
for $a \in A$ and  $b \in B$, define
$$\gamma(a,b) = \{ \gamma(a,b)_1, \gamma(a,b)_2, \ldots, 
   \gamma(a,b)_K \}$$
where 
$\gamma(a,b)_k$ equals 1 (agreement) or 0 (disagreement) on field
$k$, $k=1,\ldots,K$.  Heuristically speaking, many agreements
($\gamma(a,b)$ mostly 1's) are typical of matches, whereas many
disagreements ($\gamma(a,b)$ mostly 0's) are typical of nonmatches.  
Some variables (e.g., race) are informative in
some locations regarding matches and nonmatches, but not in others.
Disagreement on sex suggests a nonmatch, whereas agreement on sex is
not persuasive by itself for being a match.     

In this type of record linkage no one variable conclusively determines
if a pair is a match or non match.  Rather it is the composite
evidence that must be judged.  

In census and other operations, the files are divided geographically
into groups of records or 'blocks' that do not overlap.  Blocking
is used in  other applications as well in order to reduce the number
of record pairs being compared.    It is assumed that there are no (or
very few) matches across different blocks.
Other operations use first letter of last name (individuals) or
industry code (businesses) or state as blocking variables. 

Let blocks be indexed by $s=1,\ldots,S$.  Suppose that file $A$ has
$n_{a_s}$ records and file $B$ has $n_{b_s}$ records, respectively, in
block $s$.  For blocks $s=1,\ldots,S$, $a_s = 1, \ldots, n_{a_s}$ and 
$b_s = 1, \ldots, n_{b_s}$, define
$$I(a_s, b_s) = \left\{ 
  \begin{array}{cr} 1 & \hspace{1in}  a \mbox{ and } b \mbox{ are matches} \\
                    0 & a \mbox{ and } b \mbox{ are nonmatches}
  \end{array} \right.  $$
The set of match-nonmatch indicators in block $s$ is 
$I_s = \left\{ I(a_s,b_s) \right\}$.  

The match/nonmatch indicators ${\bf I} = \{ I(a,b), a \in A_s, b \in
B_s, s=1, \ldots, S \}$ are unobserved unless clerical review or some
verification system is used.  

\subsection{Latent Class Models}

The mixture model (McLachlan and Peel 2000) approach to record linkage
models the probability of 
a comparison vector $\gamma$ as arising from a mixture distribution: 
\begin{eqnarray} 
\mbox{Pr}(\gamma) &= &\mbox{Pr}(\gamma|M)p_M + \mbox{Pr}(\gamma|U)p_U,
  \label{eq:mix}
\end{eqnarray}
where $\mbox{Pr}(\gamma|M)$ and $\mbox{Pr}(\gamma|U)$ are the
probabilities of the  pattern $\gamma$ among the matches ($M$) and
nonmatches ($U$), respectively, and 
$p_M$ and $p_U = 1 - p_M$ are marginal probabilities of matches and
unmatched pairs.  

The conditional independence assumption of latent class models
(McCutcheon 1987) simplifies
the model by reducing the dimension within each mixture class from
$2^K-1$ parameters to $K$:  
\begin{eqnarray}
 \mbox{Pr}(\gamma|C) &= &\prod_{k=1}^K \mbox{Pr}(\gamma_k|C)^{\gamma_k} 
        (1-\mbox{Pr}(\gamma_k|C))^{1-\gamma_k},
  \label{eq:CI}
\end{eqnarray} 
with $C \in \{ M, U \}$.  Interactions between comparison fields have
been allowed in Larsen and Rubin (2001), Armstrong and Mayda (1993),
Thibaudeau (1993), Winkler (1989), and others.  Here we consider only the
conditional independence model and extensions of it to a hierarchical
framework. The CI assumption reduces the number of parameters needed
to describe $\mbox{Pr}(\gamma)$  from $2^K-1$ to $2K+1$.  Maximum
likelihood estimation  of parameters is accomplished with the
EM-algorithm (Dempster, Laird, and Rubin 1977).  

By Bayes' theorem, if the parameters
were known and one does not consider restrictions from one-to-one
matching, one could calculate for a pair $(a,b)$ the probability that
$a$ and $b$ match: 
\begin{eqnarray}
 \label{eq:prob1}
\mbox{Pr}(I(a,b)=1|\gamma(a,b)) \;  = 
   \; \mbox{Pr}(M| \gamma(a,b))  \; & = & \; p_M
\mbox{Pr}(\gamma(a,b)|M) / \mbox{Pr}(\gamma(a,b))
\end{eqnarray}
with the denominator given by (\ref{eq:mix}).

\subsection{Prior Information and Bayesian Latent Class Model}
Experience from previous record linkage
operations has been used informally to select models (Larsen and Rubin
2001) and restrict parameters (Winkler 1989, 1994).   Bayesian approaches to
record linkage have been suggested by Larsen (1999a, 2002, 2004,
2005), Fortini {\em et al.} (2002, 2000), and McGlinchy (2004). 

Prior experience and data often are available from previous record
linkage operations and sites.     In previous record linkage studies,
clerks at the U.S.\ Census Bureau looked at record pairs and
determined whether or not they truly were nonmatches or matches.  
Belin (1993), Belin and Rubin (1995), Larsen (1999b),  and Larsen and
Rubin (2001) found 
that in some U.S.\ Census Bureau record linkage applications
characteristics of populations being studied varied by area in ways
that made a significant impact on estimates of parameters needed for
record linkage.  There were, however, consistent patterns across
areas.  The percentage of record pairs, one record from each of two
files,  under consideration that actually are matches corresponding to
the same person is roughly similar across sites.  The probability of
agreeing on some key fields of information among matches and
nonmatches are similar across sites.  The probability of agreements
are higher among matches than among nonmatches.  There is, however,
variability across sites in these and many other characteristics.

Assuming the conditional independence model (\ref{eq:CI}) and global
parameters that do not vary by block, a prior distribution on
parameters can be  specified conveniently as the product of independent
Beta distributions as follows:  
  $$p_M \sim  \mbox{Beta}(\alpha_M, \beta_M),$$
$$\mbox{Pr}(\gamma_k(a,b)=1|M) \sim 
  \mbox{Beta}(\alpha_{Mk}, \beta_{Mk}), k=1,\ldots,K,$$ and
$$\mbox{Pr}(\gamma_k(a,b)=1|U) \sim 
  \mbox{Beta}(\alpha_{Uk}, \beta_{Uk}), k=1,\ldots,K.$$ 

Instead of specifying the prior distribution in this manner, it would
conceptually be possible to specify a prior distribution on the whole
of the probability vector associated with the set of comparison
vectors $\gamma$ as two Dirichlet distributions. That is, independent
prior distributions $\mbox{Pr}(\gamma|M) \sim
\mbox{Dirichlet}(\delta_M)$ and $\mbox{Pr}(\gamma|U) \sim
\mbox{Dirichlet}(\delta_U)$ could be specified.  This option is not
explored in this paper.  It is noted, however, that pairs of records
with known match status could be used as ``training data'' (as in
Belin and Rubin 1995) for the purposes of specifying a prior
distribution.   The prior parameter values, $\delta_M$ and $\delta_U$,
could be considered as 'prior counts' by agreement vector pattern in
the matches and nonmatches. 

If the match indicators ${\bf I}$ were known, the posterior
distributions of individual parameters given values of the other
parameters would be as follows: 
\begin{eqnarray}
 \label{eq:pM}
p_M | {\bf I} &\sim& \mbox{Beta}(\alpha_M+ 
        \sum_{(a,b)} I(a,b), \; \beta_M+ \sum_{(a,b)} (1-I(a,b)),
\end{eqnarray}
\begin{eqnarray}
 \label{eq:pkM}
\mbox{Pr}(\gamma_k(a,b)=1|M,{\bf I}) & \sim  &
 \mbox{Beta}(\alpha_{Mk} + \sum I_{ab} \gamma_k(a,b),  
 \nonumber \\
&  & \mbox{\hspace{0.5in}} \beta_{Mk}+ \sum I_{ab} (1-\gamma_k(a,b)))
\end{eqnarray}
for $k=1,\ldots,K$, 
and
\begin{eqnarray}
 \label{eq:pkU}
\mbox{Pr}(\gamma_k(a,b)=1|U,{\bf I}) & \sim &
 \mbox{Beta}(\alpha_{Uk} + \sum (1-I_{ab}) \gamma_k(a,b), 
 \nonumber   \\ 
 & & \mbox{\hspace{0.5in}} \beta_{Uk}+ \sum (1-I_{ab}) (1-\gamma_k(a,b))) 
\end{eqnarray}
for $k=1,\ldots,K$, 
where $I_{ab}=I(a,b)$ and sums are over all pairs allowed within the
blocking structure. 

\subsection{Computing for the Bayesian Latent Class Model}
The posterior distribution of parameters is simulated by sampling from
alternating conditional distributions (Gibbs sampling; Geman and Geman
1984, Gelfand and Smith 1990) as follows.  
\begin{enumerate}
\item Specify parameters for the prior distributions.  Choose initial
  values of unknown parameters. 

\item Repeat the following four steps numerous times until the distribution
of draws has converged to the posterior distribution of interest. 

\begin{enumerate}
\item Draw values for the components of  $I$ independently from
  Bernoulli distributions with the probability that $I(a,b)=1$ given
  by formula (\ref{eq:prob1}). 

\item Draw a value of $p_M$ from the distribution specified in
  formula (\ref{eq:pM}) and calculate $p_U = 1-p_M$. 

\item Draw values of $\mbox{Pr}(\gamma_k(a,b)=1|M,{\bf I})$
  independently for $k=1,\ldots,K$ from distributions specified in
  formula (\ref{eq:pkM}). 

\item Draw values of $\mbox{Pr}(\gamma_k(a,b)=1|U,{\bf I})$
  independently for $k=1,\ldots,K$ from distributions specified in
 formula (\ref{eq:pkU}). 
\end{enumerate}

\item Stop once the algorithm has converged.   Convergence of the
  algorithm can be monitored by comparing 
distributions from multiple independent series as suggested by Gelman
and Rubin (1992) and Brooks and Gelman (1998).  
\end{enumerate}

Once the algorithm has converged, it is necessary to decide which
pairs of records  to designate links and nonlinks and which to send to
clerical review or leave undecided.  
One can calculate the proportion of times that a
record pair $(a,b)$ has $I(a,b)=1$.  For record pairs with a
proportion exceeding a cut off, such as 0.90, one can make the
assignment of the pair to the match group.  Larsen (2012) examined the
impact of cutoff values. 

\subsection{Comments on Bayesian Latent Class Model}
There are some restrictions on parameters that potentially could improve the
performance of this model for record linkage.  First, the range of
$p_M$ logically should be restricted to be less than or equal to the
smaller of the two file sizes divided by the number of pairs under the
blocking structure.  When $p_M$ is drawn in the Gibbs sampling
algorithm from its conditional distribution, values of $p_M$ greater
than the cutoff should not be used. Alternatively, if $p_M = c_M p_M'$
where $p_M'$ has the Beta distribution given above and $c_M<1$
is a scale factor appropriate for transforming $p_M'$ to the allowable
range of $p_M$, one can sample $p_M'$ and scale it by $c_M$.  Second,
logically the probability of a record pair agreeing on a comparison field
should be larger among matches than among nonmatches.  That is,
$\mbox{Pr}(\gamma_k|M) > \mbox{Pr}(\gamma_k|U)$, for $k=1,\dots,K$.
Such a restriction can be added to the Gibbs sampling algorithm by
simply ignoring sampled pairs of these probabilities that do not
satisfy the constraint.   Alternatively, one can draw one value, say 
$\mbox{Pr}(\gamma_k|M)$, and scale the value of 
$\mbox{Pr}(\gamma_k|U)$ to be in the range $(0,
\mbox{Pr}(\gamma_k|M))$.  That is, after drawing a value of 
$\mbox{Pr}(\gamma_k|M)$, draw a value from the Beta distribution
specified in the algorithm and multiply it by 
$\mbox{Pr}(\gamma_k|M)$. 

Instead of specifying a  specific Beta prior distributions for
latent class parameters, one could consider specifying a hyperprior
distribution, $p(\alpha, \beta)$, for the parameters of the Beta
distributions, where $\alpha = (\alpha_M, \alpha_{Mk}, \alpha_{Uk},
k=1,\ldots,K)$ and  $\beta = (\beta_M, \beta_{Mk}, \beta_{Uk},
k=1,\ldots,K)$.  One could first transform the parameters to a scale
reflecting  proportions and sample sizes:  
$\theta = \mbox{logit}(\frac{\alpha}{\alpha+\beta})$ and $\tau =
\log(\alpha+\beta)$, where the transformation is applied to
corresponding components of the $\alpha$ and $\beta$ vectors. 
Note that there is a unique bivariate inverse transformation:
$\alpha = e^{\tau} \mbox{logit}^{-1}(\theta)$
and 
$\beta = e^{\tau} \mbox{logit}^{-1}(1-\theta)$. 
The approach of specifying a hyperprior distribution is not considered
here, because such an approach implicitly makes the assumption that
the parameters are exchangeable.  Such an assumption is unrealistic
because variables have very different levels of agreement for both
matches and nonmatches.

\section{Record Linkage and Hierarchical Latent Class Model}

A hierarchical model for record linkage specifies distributions of
parameters within blocks $s=1,\ldots,S$.  
The probabilities of agreeing on fields of information are allowed to
vary by block.   Prior distributions on parameters are as follows: 
$$p_{sMk} = \mbox{Pr}(\gamma_k=1| M, s)  \sim
\mbox{Beta}(\alpha_{sMk}, \beta_{sMk})$$ 
and 
$$ p_{sUk} = \mbox{Pr}(\gamma_k=1| U, s)  \sim
\mbox{Beta}(\alpha_{sUk},\beta_{sUk})$$
independently across blocks ($s=1,\ldots,S$), fields ($k=1,\ldots,K$),
and classes ($M$ and $U$).  
As before one can assume that the restriction that $p_{sMk} \ge
p_{sUk}$. 

\subsection{Hyperprior distributions for the Hierarchical Model}
Hyperprior distributions are used in this model to link estimation of
probabilities across blocks.  Without a model that enables 'borrowing
strength' across blocks, parameters appear separately by blocks and
data could be insufficient for accurately estimation.  It is still
the case that it is unlikely to be reasonable to model
the probability of matching and the probability of
agreement for $K$ fields as exchangeable.   As a result, independent
hyperprior distributions are used.  One version was suggested by
Larsen (2004).  The specification below generalizes previous
approaches by allowing correlation.  Let 
 $$\theta_{sMk} = 
 \mbox{logit}(\frac{\alpha_{sMk}}{\alpha_{sMk}+\beta_{sMk}}),$$ 
 $$\tau_{sMk}  = \log (\alpha_{sMk} +  \beta_{sMk}),$$ 
 $$ \theta_{sUk} =
  \mbox{logit}(\frac{\alpha_{sUk}}{\alpha_{sUk}+\beta_{sUk}}),$$ 
 $$\tau_{sUk} = \log (\alpha_{sUk} +  \beta_{sUk}),$$ 
 $$\theta_{sM} =
   \mbox{logit}(\frac{\alpha_{sM}}{\alpha_{sM}+\beta_{sM}}),$$ 
and 
 $$\tau_{sM} = \log (\alpha_{sM} +  \beta_{sM}). $$   

Then 
 $$ (\theta_{sMk}, \tau_{sMk})^{t} \sim 
    N_2(  (\mu_{\theta M k}, \mu_{\tau M k})^{t}, \Sigma_{Mk} ),$$

 $$ (\theta_{sUk}, \tau_{sUk})^{t} \sim 
    N_2(  (\mu_{\theta U k}, \mu_{\tau U k})^{t}, \Sigma_{Uk} ),$$
and
 $$ (\theta_{sM}, \tau_{sM})^{t} \sim 
    N_2(  (\mu_{\theta M}, \mu_{\tau M})^{t}, \Sigma_{M} ),$$
indepedendently,  where 
$$ \Sigma_{Mk} = \left( 
\begin{array}{cc}
\sigma^2_{\theta M k} & \sigma_{\theta M k, \tau M k} \\
\sigma_{\theta M k, \tau M k}  & \sigma^2_{\tau M k}
\end{array} 
 \right), $$ 
$$ \Sigma_{Uk} = \left( 
\begin{array}{cc}
\sigma^2_{\theta U k} & \sigma_{\theta U k, \tau U k} \\
\sigma_{\theta U k, \tau U k}  & \sigma^2_{\tau U k}
\end{array} 
 \right), $$ 
and
$$ \Sigma_{M} = \left( 
\begin{array}{cc}
\sigma^2_{\theta M} & \sigma_{\theta M, \tau M} \\
\sigma_{\theta M, \tau M}  & \sigma^2_{\tau M}
\end{array} 
 \right). $$ 

In the prior distributions, one could enforce the restriction 
that, for $k=1,\ldots,K$, $ \theta_{sMk} \ge \theta_{sUk}$. 
Similarly the restriction that
$p_{sM}$ is smaller than the minimum of $n_{A_s}$ and $n_{B_s}$
divided by the number of pairs $n_{A_s} n_{B_s}$ likely will be
useful.  If it were not enforced, the small sample size and great
variability across blocks would surely produce poor results for some
blocks.

 \subsection{Computing for the Hierarchical Latent Class Model    \label{sect:comp}}

The posterior distribution of parameters and unobserved match/nonmatch
indicators will be simulated using Gibbs sampling.  The conditional
distributions for the hyperparameters will be sampled using the
Metropolis-Hastings (MH) algorithm (Hastings 1970) within the Gibbs
sampling framework.  
The procedure iterates through draws of full conditional distributions
as described below. 

\begin{enumerate}
\item Choose hyperparameter distributions.  That is, select values for
  means ($\mu_{\theta M}$, $\mu_{\theta M k}$, $(\mu_{\theta U k}$, 
        $\mu_{\tau M}$, $\mu_{\tau M k}$, and $\mu_{\tau U k}$)  and
        variance matrices ($\Sigma_{M}$, $\Sigma_{Mk}$, $\Sigma_{Uk}$).

\item Generate initial values of  ($\alpha_{sM}, \beta_{sM}$)
and, for $k=1,\ldots,K$,
$(\alpha_{sMk}, \beta_{sMk})$ and 
$(\alpha_{sUk}, \beta_{sUk})$
from their prior distributions. 

\item Assign an initial match/nonmatch configuration ${\bf I}$.   One
  initialization option is to 
  randomly generate a matrix of 1's and 0's representing matches and
  nonmatches by block with the number of
  1's per block below the maximum allowable number of 1's.  A possibly
  better way to initialize ${\bf I}$ would be to use the ordinary
  latent class model to assign matches and nonmatches at a high
  probability cutoff. 

\item Cycle through the following steps numerous times until the
  distribution of drawn values converges to the target posterior
  distribution.   Let $I_{ab}$ denote $I(a,b)$. 

 \begin{enumerate}
 \item  For $s=1,\ldots,S$, draw   $p_{sM}$ from its conditional
   distribution given the current indicators ${\bf I}_s$
   and values of
   $(\alpha_{sM}, \beta_{sM})$.   Specifically,
   $$p_{sM}|I_s,  \alpha_{sM}, \beta_{sM} \sim 
   \mbox{Beta}(\alpha_{sM}+ \sum I_{ab}, \beta_{sM} + n_{a_s} n_{b_s} -
   \sum I_{ab}), $$
    where the sum is over all pairs $(a,b)$ in block $s$.   Enforce
    the constraint $$p_{sM} \le \mbox{min}(n_{a_s}, n_{b_s})/ 
    (n_{a_s} n_{b_s}).$$

 \item For $s=1,\ldots,S$ and $k=1,\ldots,K$,  draw 
   $p_{sMk}$ and $p_{sUk}$ from their conditional
   distribution given the current indicators ${\bf I}_s$, the
   comparison vectors $\gamma_s$ in block $s$,  and values of
   $(\alpha_{sCk}, \beta_{sCk}), C \in \{M,U\}$.   Specifically,
   $$p_{sMk}|I_s, \gamma_s, \alpha_{sMk}, \beta_{sMk} \sim 
   \mbox{Beta}(\alpha_{sMk}+ \sum_s I_{ab} \gamma_k(a,b), 
             \beta_{sMk} + \sum_s I_{ab} (1-\gamma_k(a,b))),$$
   $$p_{sUk}|I_s, \gamma_s, \alpha_{sUk}, \beta_{sUk} \sim 
   \mbox{Beta}(\alpha_{sUk}+ \sum_s (1-I_{ab}) \gamma_k(a,b), 
             \beta_{sUk} + \sum_s (1-I_{ab}) (1-\gamma_k(a,b))),$$
  and 
    $p_{sMk} \ge p_{sUk}$, 
    where sums are over all pairs $(a,b)$ in block $s$. 

 \item  For $s=1,\ldots,S$, use the Metropolis-Hastings algorithm
   (Hastings 1970; see also 
   Gelman 1992 and Gelman {\em et al.} 2004, chapter 11) to draw values of
   hyperparameters $\theta_{sM}$ and $\tau_{sM}$ from their full
   conditional distributions.   Details of this step and the next two
   steps are given after this outline. 

\item For $s=1,\ldots,S$ and $k=1,\ldots,K$, use the
  Metropolis-Hastings algorithm  to 
  draw values of hyperparameters $\theta_{sMk}$ and $\tau_{sMk}$.

\item For $s=1,\ldots,S$ and $k=1,\ldots,K$, use the
  Metropolis-Hastings algorithm  to 
  draw values of hyperparameters $\theta_{sUk}$ and $\tau_{sUk}$.

 \item For $s=1,\ldots,S$, $a=1,\ldots,n_{a_s}$, and
   $b=1,\ldots,n_{b_s}$, given values of $p_{sM}$ and, for
   $k=1,\ldots,K$, $p_{sMk}$ and $p_{sUk}$, draw a value of $I(a,b)$
   from a Bernoulli distribution with the following probability: 
$$
\frac{ p_{sM} \prod_{k=1}^K \left[ p_{sMk}^{\gamma_k(a,b)}
        (1-p_{sMk})^{1-\gamma_k(a,b)} \right]   }
{  \left\{ p_{sM} \prod_{k=1}^K \left[ p_{sMk}^{\gamma_k(a,b)}
        (1-p_{sMk})^{1-\gamma_k(a,b)} \right] + 
(1-p_{sM}) \prod_{k=1}^K \left[ p_{sUk}^{\gamma_k(a,b)}
        (1-p_{sUk})^{1-\gamma_k(a,b)} \right] \right\}  }.$$

 \end{enumerate}

\item Stop once the algorithm has converged. 
\end{enumerate}

As before, once the algorithm has converged, it is necessary to decide which
pairs of records  to designate as links.  Suggestions were made at the end
of the previous section. 

Details of the three Metropolis-Hastings (Hastings 1970) steps in the
simulation procedure above are now presented. 

\renewcommand{\theenumi}{(\alph{enumi})}

\begin{enumerate}
\setcounter{enumi}{2}
 \item  \label{appAc} For $s=1,\ldots,S$, use the Metropolis-Hastings algorithm
    to draw values of
   hyperparameters $\theta_{sM}$ and $\tau_{sM}$ from their full
   conditional distributions.  Specifically, given current values of
   $\theta_{sM}$ and  $\tau_{sM}$ (and hence $\alpha_{sM}$ and
   $\beta_{sM}$), ${\bf I}_s$, and other parameters,  

 \renewcommand{\theenumii}{\roman{enumii}}
 \begin{enumerate}
 \item Define a tuning constant $h_{M}>0$. 
 \item Draw three values: 
  $$u \sim \mbox{Uniform}(0,1)$$
and 
  $$(\theta^*, \tau^*)^{t} \sim N_2(
    (\theta_{sM}, \tau_{sM})^{t},  \Sigma_M/h_M).$$

\item Calculate $\alpha^* =  e^{\tau^*}
  \mbox{logit}^{-1}(\theta^*)$ and 
  $\beta^* =  e^{\tau^*}
  \mbox{logit}^{-1}(1-\theta^*)$. 

\item Calculate 
\begin{eqnarray*}
r & = & \mbox{min} \left(
 1, \;   p_{sM}^{\alpha^* - \alpha_{sM}} 
       (1-p_{sM})^{\beta^* - \beta_{sM}} \right.  \nonumber \\
  & &  \times   \left. \exp{( -\frac{h_{\theta M}}{\sigma^2_{\theta M}} 
                 (\theta_{sM}-\theta^*)^2 )}
       \exp{( -\frac{h_{\tau M}}{\sigma^2_{\tau M}} 
                 (\tau_{sM}-\tau^*)^2 )}
 \right) 
\end{eqnarray*}

\item 
If $u \le r$, let $\theta_{sM} = \theta^*$ and 
 $\tau_{sM} = \tau^*$.  

Otherwise, let $\theta_{sM}$ and $\tau_{sM}$ remain the same. 

\end{enumerate}

\item For $s=1,\ldots,S$ and $k=1,\ldots,K$, use the
  Metropolis-Hastings algorithm  to 
  draw values of hyperparameters $\theta_{sMk}$ and $\tau_{sMk}$. 
Specifically, given current values of
   $\theta_{sMk}$ and  $\tau_{sMk}$ (and hence $\alpha_{sMk}$ and
   $\beta_{sMk}$), ${\bf I}_s$, and other parameters,  
  follow the   steps outlined in step \ref{appAc} above 
  but with all $M$ indexes    replaced by $Mk$'s.  The tuning
  parameter $h_{Mk}>0$ needs to be chosen. 







\item For $s=1,\ldots,S$ and $k=1,\ldots,K$, use the
  Metropolis-Hastings algorithm  to 
  draw values of hyperparameters $\theta_{sUk}$ and $\tau_{sUk}$.
Specifically, given current values of
   $\theta_{sUk}$ and  $\tau_{sUk}$ (and hence $\alpha_{sUk}$ and
   $\beta_{sUk}$), ${\bf I}_s$, and other parameters,  follow the
   steps outlined in step \ref{appAc} above 
  but with all $M$  indexes    replaced by $Uk$'s. 
    The tuning   parameter $h_{Uk}>0$ needs to be chosen. 







\end{enumerate}



The tuning parameters $h_{M}$ and, for $k=1,\ldots,K$, $h_{Mk}$ and
$h_{Uk}$ are chosen so 
that the drawn values of the parameters are accepted approximately
about 35\% of the time (Gelman {\em et al.} 2004 chapter 11.9).  The
algorithm could be run for several iterations to assess the acceptance
rate, adapting the tuning parameters as necessary.  A second phase then
could be initiated with fixed values for tuning parameters.  

\section{Simulation \label{sect:sim}}

One thousand replications are performed under each of two sets of
conditions.  In one set of conditions, the probabilities of agreement
are constant across blocks.  In the second set of conditions, the
probabilities of agreement vary by block.  Blocks are assumed to be
linked together correctly, as they would be if they correspond to
geographical areas.  Pairs from different blocks are nonlinks and not
used to estimate probabilities. Files $A$ and $B$ both have 10,000
records organized into 400 blocks of size 25 each.   This arrangement
yields 250,000 (400 times $25^2$) pairs of records.  

In the first scenario, probabilities of matching and not matching are
constant across blocks. The seven matching variables have probabilities of
agreement among matches of 0.90 to 0.60 in increments of 0.05.   The
probabilities of agreement among nonmatches is 0.5 for two variables,
0.33 for three variables, and 0.25 for two variables.  Agreements
on the fields of information are independent  of one another.

In the second scenario, probabilities of matching and not matching 
vary across blocks. The seven matching variables in each block have
probabilities randomly selected from the range 0.60 to 0.90.    The
probabilities of agreement among nonmatches have probabilities
randomly selected from 0.20 to 0.50. As before, agreements
on the fields of information are independent  of one another.

The files $A$ and $B$ were generated and comparison vectors
calculated.  

Four statistical procedures were used to do the record linkage.  The
first is latent class analysis (LCA).  The second is the Bayesian
latent class model (BLCM) with a chosen prior distribution.   The third is
the BLCM with a prior distribution based on a sample of similar cases.
The fourth is the hierarchical LCA. 

For the LCA, the EM algorithm (Dempster, Laird, Rubin 1977) was used to fit a
two-class conditional independence mixture model to the comparison
vectors to estimate probabilities for the Fellegi-Sunter (1969)
algorithm. 

In the first Bayesian LCM, Beta prior distributions are chosen so that the
probability of agreement among matches, $\mbox{Pr}(\gamma_k=1|M)$, is
most likely between 0.65 and 0.95, the probability among nonmathces,
$\mbox{Pr}(\gamma_k=1|U)$, tends to be between 0.10 and 0.40, and the
probability of a matching pair, $\mbox{Pr}(M)$, is likely
between 0.02 and 0.04.   If the mean of a Beta
distribution is 0.80 and its standard deviation (SD) is 0.075, then its
parameters are approximately $\alpha_{Mk} = 22.0$ and $\beta_{Mk} =
5.4$.   A Beta distribution with mean 0.25 and SD
0.075 has parameters approximately $\alpha_{Uk} = 8.1$ and $\beta_{Mk} =
42.2$.   A Beta distribution with mean 0.03 and SD 0.005 has
parameters approximately $\alpha_{M} = 35$ and $\beta_{M} = 1128$.
The prior distribution for $\mbox{Pr}(M)$ is small narrow because
logically the percent of matches has to be below 4\%:  10,000/250,000
= 1/25 = 0.04.  

In the second Bayesian LCM, it is assumed that there are complete data
from two blocks.  That is, match status is known for all 1,250 pairs
in the two blocks.  Results of these comparisons are used as prior
information.  The number of matches and nonmathing pairs determine
$\alpha_M$ and $\beta_M$, respectively.  The number of agreements and
disagreements on field $k (k=1,\ldots,K)$ among matches produces
$\alpha_{Mk}$ and $\beta_{Mk}$, respectively.  Similarly,  number of
agreements and disagreements on field $k (k=1,\ldots,K)$ among non
matches produces $\alpha_{Uk}$ and $\beta_{Uk}$, respectively.  The
sum of $\alpha_M$ and $\beta_M$ will be 1,250, which is comparable to
their sum in the first BLCM approach. 

In the hierarchical LCA, it is necessary to select hyperprior
distributions and tuning constants for the Metropolis-Hastings
steps. If one uses the means of the prior distributions from the first
Bayesian LCM formulation for the hyperprior means, then
$\mu_{\theta M} = \mbox{logit}(0.80) = 1.39$, 
$\mu_{\theta M k} = \mbox{logit}(0.25) = -1.10$, 
$(\mu_{\theta U k} = \mbox{logit}(0.03) = -3.48$ , 
$\mu_{\tau M} = \log{27.4} = 3.3$, 
$\mu_{\tau M k} = \log{32.3} = 3.5$, and 
$\mu_{\tau U k} = \log{1163} = 7.1$.  

Hyperprior variances and covariances describe spread and correlation
of $\alpha$ and $\beta$ values (transformed to $\theta$ and $\tau$)
across blocks. An ad hoc way of specifying these values is suggested.
Among matches under the prior distribution from the 
first BLCM, one expects probabilities to be between 0.65 and 0.90.  On
the logit scale, 0.65 is 0.61 and 0.90 is 2.20.  A uniform
distribution bewteen 0.61 and 2.20 has variance (2.20-0.61)/12 =
0.1325.  The prior ``sample size'' from the first BLCM was 27.4.  One
quarter of this is 6.85, which is 1.92 on the log scale.  Four times
this value is 109.60, which is 4.70 on the log scale. A uniform
distriution between 1.92 and 4.70 has variance (4.70-1.92)/12 = 0.23.
Thus $\sigma^2_{\theta Mk} = 0.1325$ and $\sigma^2_{\tau Mk} = 0.23$
for $k=1,\ldots,K$.  

Similar methods can be used for the other variances.    The logit of
0.10 is -2.20.  The logit of 0.40 is -0.41.  The variance of a uniform
distribution with these limits is 0.1492.   One quarter of 32.3 is 8.1,
which is 2.1 on the log scale.  Four times 32.3 is 129.2, which is
4.86 on the log scale.  The variance of a uniform with these limits is
0.23.   Thus  $\sigma^2_{\theta Uk} = 0.1492$ and $\sigma^2_{\tau Uk} = 0.23$
for $k=1,\ldots,K$.      The logits of 0.02 and 0.04 are -3.89 and
-3.18, respectively.   One-twelfth that range is 0.059.  One quarter
of 1163 is 291, which is 5.67 on the log scale.  Four times it is
4652, which is 8.45 after log transformation.  One-twelfth the range
is 0.23.  Thus $\sigma^2_{\theta M} = 0.059$ and $\sigma^2_{\tau M} =
0.23$. 

For the covariance terms between $\theta$ and $\tau$ parameters, an ad
hoc process was followed.   Values of $\alpha$ and $\beta$ were
computed for a range of Beta distributions with probabilities in
various ranges.  These values were transformed to $\theta$ and
$\beta$.  The covariance among them was computed.  This was repeated
for the match probability and probabilities of agreement among matches
and among nonmatches.  The values obtained are as follows:
$\sigma_{\theta Mk, \tau Mk} = -0.08$, 
$\sigma_{\theta Uk, \tau Uk} = -0.01$,  and
$\sigma_{\theta M, \tau M} = 0.03$. 

Choices for the tuning constants ($h_M$ and for $k=1,\ldots,K$
$h_{Mk}$ and $h_{Uk}$) also need to be made.  Initially these will be
set of 0.5.  Depending on acceptance rates in the initial set of
Metropolis-Hastings steps, these values will be reassessed.  The
algorithm could use $2K+1$ different values to improve algorithm
performance.

Work on simulations is underway and should be completed in early
2013. 

\section{Conclusions and Future Work \label{sect:concl}}
A novel hierarchical Bayesian model for record linkage has been
presented.  The model allows probabilities to vary by
block and reflect local information.  
Simulations are being completed to evaluate the performance of the proposed
methods. 

Several areas can be identified for future work.  Many of these will
be important in actual applications.  It will be interesting to apply
these methods to data from the U.S.\ Census Bureau, the U.S.\
National Center for Health Statistics,  and  
other sources.  An automated system for applying these models to new
sets of files would be useful in this regard.   In a real application,
one could consider better specifications of prior distributions for
the record linkage model parameters.   In particular, if data are
available from another record linkage site and the site differs in
some ways from the current application, then one must decide the
degree to which data from the previous site should be discounted or down
weighted when analyzing the new site.   In some applications, the size
of the files will be a challenge.  In order to speed computations, one
might consider parallel computations; for example, many computations
are performed separately in each block.  

The algorithm's performance could be improved by studying tuning
parameters and the order of sampling cycles within Metropolis-Hastings
and Gibbs sampling  algorithms.   One could study the sensitivity of
results to the specification of hyperprior distributions.  

Larsen (2004, 2005) considered one-to-one restrictions enforced in the
indicator matrix ${\bf I}$ and in the statistical likelihood function.
 The one-to-one restrictions and blocking assumptions mean that
$\sum_{b_s} I(a_s,b_s) \le 1$, $\sum_{a_s} I(a_s,b_s) \le 1$, and
$\sum_{a_s} \sum_{b_{s'}} I(a_s, b_{s'}) = 0$ for $s \ne s'$.  The
number of matches in block $s$, $n_{m_s}$ is defined and restricted
under one-to-one matching as
follows: 
$$\sum_{a_s} \sum_{b_s} I(a_s,b_s) = n_{m_s} \le
\min{(n_{a_s},n_{b_s})}.$$ 
Future work will pursue Metropolis-Hastings steps for sampling new
values of ${\bf I}$ instead of the current simpler formulation.  

\section*{Acknowledgments}
The work contained in this report was partially supported by NSA Grant
H98230-10-1-0232 and by a Selective Excellence Award from The George
Washington University.   The opinions expressed are those of the
author and not necessarily of the National Security Agency, The George
Washington University, or anyone else.  The author wishes to thank
Yuan Zhao and Jing Zhang for discussions on record linkage.  The
author also thanks the Department of Biostatistics of the University
of North Carolina at Chapel Hill, Westat, the Federal Committee on
Statistical Methodology, the U.S.\ National Center for Health
Statistics, and the U.S.\ Census Bureau for opportunities to present
seminars and discuss record linkage. 

\section*{References \label{s.ref}}
\noi Alvey, W., and Jamerson, B. (1997), {\em Record Linkage
  Techniques -- 1997},  Proceedings of an International Workshop and
  Exposition. Federal Committee on Statistical Methodology, Office of
  Management of the Budget. \\

\noi Armstrong, J.\ B., and Mayda, J.\ E. (1993).
Model-Based Estimation of Record Linkage Error Rates. 
{\em Survey Methodology}, 19, 137-147. \\%

\noindent
Belin, T.\ R. (1993). 
Evaluation of sources of variation in record linkage through a
factorial experiment. {\em Survey Methodology} {\bf 19}, 13-29. 
\newline


\noindent 
Belin, T.\ R., and Rubin, D.\ B. (1995). 
A method for calibrating false match rates in record linkage.
{\em Journal of the American Statistical Association} {\bf 90}, 694- 707.
\newline

\noindent 
Brooks, S.\ P. , and Gelman, A.  (1998). General methods for
monitoring convergence of iterative simulations.  {\em Journal of
Computational and Graphical Statistics}, 7 , 434-455   \\

\noi Dempster, A.\ P., Laird, N.\ M., and Rubin, D.\ B. (1977),
``Maximum Likelihood from Incomplete Data Via the EM Algorithm, ''
(with comments),
{\em Journal of the Royal Statistical Society}, Ser.\ B,
39, 1-37. \\%

\noi Fellegi, I.\ P., and Sunter, A.\ B. (1969), ``A Theory for Record
Linkage,'' {\em   Journal of the American Statistical Association},
64, 1183-1210. \\

\noindent 
Fortini, M.,  Liseo, B., Nuccitelli, A., and Scanu, M. (2000).
``On Bayesian Record Linkage,'' {\em Bayesian Methods with
  Applications to Science, Policy, and Official Statistics}: 
Selected Papers from ISBA 2000: The Sixth World Meeting of the
International Society for Bayesian Analysis.  Editor E.\ I.\ George,
155-164. \\

\noindent 
Fortini, M.,  Nuccitelli, A., Liseo, B., and Scanu, M. (2002).
``Modelling issues in record linkage: A Bayesian perspective,'' 
 {\em Proceedings of the American Statistical Association, Survey
   Research Methods Section}, 1008-1013.  \\ 

\noi Gelfand, Alan E. , and Smith, Adrian F. M.  (1990),
``Sampling-based approaches to calculating marginal densities'',
Journal of the American Statistical Association, 85 , 398-409  \\

\noi  Gelman, A. (1992). Iterative and non-iterative simulation
algorithms. {\em   Computing Science and Statistics} {\bf 24},
433-438. \\

\noi
Gelman, A., Carlin, J.\ B., Stern, H.\ S., and Rubin, D.\
B. (2004). {\em   Bayesian Data Analysis}, 2nd edition.  Chapman \&
Hall/CRC. \\

\noindent 
Gelman, A., and Rubin, D.B. (1992). 
Inference from iterative simulation using multiple sequences.
{\em Statistical Science} 7, 457- 472.  \newline

\noi Geman, S., and Geman, D.  (1984). Stochastic
relaxation, Gibbs distributions, and the Bayesian restoration of
images.  {\em IEEE Transactions on Pattern Analysis and Machine
Intelligence}, 6 , 721-741  \\

\noi Hastings, W.\ K. (1970).  Monte Carlo sapling methods using
Markov chains and their applications, {\em Biometrika}, 57, 97-109. \\

\noi Herzog, T.N., Scheuren, F.J., and Winkler, W.E.  (2007). {\em Data
quality and record linkage techniques}. New York: Springer
Science+Business Media. \\

\noi Jaro, M.\ A. (1989), ``Advances in Record-Linkage Methodology as
Applied to Matching the 1985 Census of Tampa, Florida,'' {\em Journal
  of the American Statistical Association},  84, 414-420.\\

\noi Jaro, M.\ A. (1995), ``Probabilistic Linkage of Large Public Health
Data Files, '' {\em Statistics in Medicine}, 14, 491-498. \\

\noi Lahiri, P.,  and Larsen, M.D. (2005). Regression Analysis with
Linked Data. {\em Journal of the  American Statistical   Association}, 100,
222-230.  \\  

\noi Larsen, M.D. (1999a), ``Multiple Imputation Analysis of Records
Linkage Using Mixture Models,,'' {\em Proceedings of the Statistical
  Society of Canada, Survey Methods Section}, 65-71. \\

\noindent Larsen, M.D. (1999b). ``Predicting the Residency Status for
Administrative Records    that Do Not Match Census Records,'' 
{\em Administrative   Records Research Memorandum Series}, \#20,
Bureau of the Census, U.S.\   Department of Commerce. \newline

\noindent 
Larsen, M.D.  (2002),
``Comment on Hierarchical Bayesian Record Linkage,''
{\em Proceedings of the Section on Bayesian Statistical Science}, American
Statistical Association meeting in New York City.  CDROM:
1995-2000. \newline 

\noindent
Larsen, M.D. (2004), Record linkage using finite mixture models, 
{\em Applied Bayesian Modeling and Causal
  Inference from Incomplete-Data Perspectives}. Gelman, A., and Meng,
X.\ L., editors. New York: Wiley, 309-318. \\

\noi
Larsen, M.D. (2005), Hierarchical Bayesian record linkage.  Technical
report, Iowa State University, Department of Statistics. \\

\noi
Larsen, M.D. (2012), Factors affecting record linkage to the National
Death index using Latent Class Models.  {\em Journal of Survey
  Statistics and Methodology}, submitted. \\

\noi  Larsen, M.\ D., and Rubin, D.\ B. (2001), ``Iterative automated
record linkage using mixture models,''         {\em Journal of the American
    Statistical Association}, 96, 32-41. \\%

\noindent McCutcheon, A.\ L.\  (1987). {\em Latent class analysis}. Sage
Publications, Inc.: Newbury Park, CA; London. \\

\noi
McGlinchy, M.\ (2004). A Bayesian record linkage methodology for
multiple imputation of missing links.  
{\em Proceedings of the American Statistical Association, Section on
  Survey Research Methods}.  Alexandria, VA: CDROM.  \\

\noi
McLachlan, G.J., and Peel, D. (2000). {\em Finite Mixture Models}. New
York: Wiley. \\





\noi Thibaudeau, Y. (1993), ``The Discrimination Power of
Dependency Structures in Record Linkage," {\em Survey
Methodology}, 19, 31-38.\\

\noi Winkler, W.\ E. (1988),
``Using the EM Algorithm for Weight Computation in the
Fellegi-Sunter Model of Record Linkage,''
{\em American Statistical Association Proceedings of Survey Research
Methods Section},  pp.\ 667- 671. \\%

\noi
 Winkler, W.\ E.  (1989), ``Near automatic weight computation in
 the Fellegi-Sunter model of record linkage'', {\em Proceedings of the
 Bureau of the Census Annual Research Conference}, 5, 145-155. \\



\noi Winkler, W.\ E. (1994),
``Advanced Methods for Record Linkage, '' in
{\em American Statistical Association Proceedings of Survey Research
Methods Section},  pp.\ 467-472. \\%

\noi Winkler, W.\ E.  (1995), ``Matching and Record
Linkage,'' in {\em Business Survey Methods}, ed.\ Cox, B.\ G.,
Binder, D.\ A., Chinnappa, B.\ N., Christianson, A.,
Colledge, M.\ J., and Kott, P.\ S., New York: Wiley Publications, pp.\
355-384.

\end{document}